\documentclass[twoside,leqno,10pt]{article}
\usepackage{graphics,ifthen}
\usepackage{latexsym}
\usepackage{amsmath}
\usepackage{amssymb}
\usepackage{mathrsfs}
\begin{document}
\input{latexP.sty}
\input{referencesP.sty}
\input epsf.sty

\def\ind{\stackrel{\mathrm{ind}}{\sim}}
\def\iid{\stackrel{\mathrm{iid}}{\sim}}

\def\Definition{\stepcounter{definitionN}\
    \Demo{Definition\hskip\smallindent\thedefinitionN}}
\def\EndDefinition{\EndDemo}
\def\Example#1{\Demo{Example [{\rm #1}]}}
\def\EndExample{\qed\EndDemo}
\def\Category#1{\centerline{\Heading #1}\rm}
\
%% Paper specific definitions
\def\e{\text{\hskip1.5pt e}}
\newcommand{\eps}{\epsilon}
\newcommand{\proof}{\noindent {\bf Proof:\ }}
\newcommand{\remarks}{\noindent {\bf Remarks:\ }}
\newcommand{\note}{\noindent {\bf Note:\ }}
\newcommand{\examp}{\noindent {\bf Example:\ }}
\newcommand{\Lower}[2]{\smash{\lower #1 \hbox{#2}}}
\newcommand{\ben}{\begin{enumerate}}
\newcommand{\een}{\end{enumerate}}
\newcommand{\bi}{\begin{itemize}}
\newcommand{\ei}{\end{itemize}}
\newcommand{\hp}{\hspace{.2in}}

\newtheorem{lw}{Proposition 3.1, Lo and Weng (1989)}
\newtheorem{thm}{Theorem}[section]
\newtheorem{defin}{Definition}[section]
\newtheorem{prop}{Proposition}[section]
\newtheorem{lem}{Lemma}[section]
\newtheorem{cor}{Corollary}[section]
\newcommand{\rb}[1]{\raisebox{1.5ex}[0pt]{#1}}
\newcommand{\mc}{\multicolumn}
%Mathrsfs Font
\newcommand{\Bcr}{\mathscr{B}}
\newcommand{\Ucr}{\mathscr{U}}
\newcommand{\Gcr}{\mathscr{G}}
\newcommand{\Dcr}{\mathscr{D}}
\newcommand{\CS}{\mathscr{C}}
\newcommand{\Fcr}{\mathscr{F}}
\newcommand{\Icr}{\mathscr{I}}
\newcommand{\Lcr}{\mathscr{L}}
\newcommand{\Mcr}{\mathscr{M}}
\newcommand{\Ncr}{\mathscr{N}}
\newcommand{\Pcr}{\mathscr{P}}
\newcommand{\Qcr}{\mathscr{Q}}
\newcommand{\Scr}{\mathscr{S}}
\newcommand{\Tcr}{\mathscr{T}}
\newcommand{\Xcr}{\mathscr{X}}
\newcommand{\Vcr}{\mathscr{V}}
\newcommand{\Ycr}{\mathscr{Y}}
%Mathbb Font
\newcommand{\E}{\mathbb{E}}
\newcommand{\F}{\mathbb{F}}
\newcommand{\I}{\mathbb{I}}
\newcommand{\Q}{\mathbb{Q}}
\newcommand{\X}{\mathbb{X}}
\newcommand{\Pe}{\mathbb{P}}
\newcommand{\M}{\mathbb{M}}
\newcommand{\Wbb}{\mathbb{W}}

\def\Beta{\text{Beta}}
\def\Dir{\text{Dirichlet}}
\def\DP{\text{DP}}
\def\P{{\bf p}}
\def\fhat{\widehat{f}}
\def\GA{\text{gamma}}
\def\ind{\stackrel{\mathrm{ind}}{\sim}}
\def\iid{\stackrel{\mathrm{iid}}{\sim}}
\def\J{{\bf J}}
\def\K{{\bf K}}
\def\min{\text{min}}
\def\N{\text{N}}
\def\p{{\bf p}}
\def\U{{\bf U}}
\def\W{{\bf W}}
\def\T{{\bf T}}
\def\y{{\bf y}}
\def\t{{\bf t}}
\def\m{{\bf m}}
\def\X{{\bf X}}
\def\Y{{\bf Y}}
\def\tps{\mbox{\scriptsize ${\theta H}$}}   %   smaller "\psi"-vector
\def\ups{\mbox{\scriptsize ${P_{\theta}(g)}$}}   %   smaller "\psi"-vector
\def\vps{\mbox{\scriptsize ${\theta}$}}   %   smaller "\psi"-vector
\def\vups{\mbox{\scriptsize ${\theta >0}$}}   %   smaller "\psi"-vector
\def\hps{\mbox{\scriptsize ${H}$}}   %   smaller "\psi"-vector
\def\rps{\mbox{\scriptsize ${(\theta+1/2,\theta+1/2)}$}}   %   smaller "\psi"-vector
\def\sps{\mbox{\scriptsize ${(1/2,1/2)}$}}   %   smaller "\psi"-vector

\newcommand{\reals}{{\rm I\!R}}
\newcommand{\PR}{{\rm I\!P}}
\def\Z{{\bf Z}}
\def\yy{{\mathcal Y}}
\def\rr{{\mathcal R}}
\def\BP{\text{beta}}
\def\ts{\tilde{t}}
\def\js{\tilde{J}}
\def\gs{\tilde{g}}
\def\fs{\tilde{f}}
\def\ys{\tilde{Y}}
\def\ps{\tilde{\mathcal {P}}}

\def\Report{Lancelot F. James}
\def\Author{Integrated Volatility}
\pagestyle{myheadings}
\markboth{\Author}{\Report}
\thispagestyle{empty}

\bct\Heading  Analysis of a Class of Likelihood Based Continuous
Time Stochastic Volatility Models including Ornstein-Uhlenbeck
Models in Financial Economics.\lbk\lbk\smc Lancelot F.
James\footnote{ \eightit AMS 2000 subject classifications.
               \rm Primary 62G05; secondary 62F15.\\
\eightit Corresponding authors address.
                \rm The Hong Kong University of Science and Technology,
Department of Information and Systems Management, Clear Water Bay,
Kowloon, Hong Kong.
\rm lancelot\at ust.hk\\
\indent\eightit Keywords and phrases.
                \rm
          Bessel Functions,
          Mixture of Normals,
          Ornstein-Uhlenbeck Process,
          Poisson Process,
          Stochastic Volatility,
          Weber-Sonine Formula.
          }
\lbk\lbk \BigSlant The Hong Kong University of Science and
Technology\rm \lbk %(\today)%
\ect \Quote In a series of recent papers Barndorff-Nielsen and
Shephard introduce an attractive class of continuous time
stochastic volatility models for financial assets where the
volatility processes are functions of positive
Ornstein-Uhlenbeck(OU) processes. This models are known to be
substantially more flexible than Gaussian based models. One
current problem of this approach is the unavailability of a
tractable exact analysis of likelihood based stochastic volatility
models for the returns of log prices of stocks.
 With this point in mind, the likelihood models of
Barndorff-Nielsen and Shephard are viewed as members of a much
larger class of models. That is likelihoods based on n
conditionally independent Normal random variables whose mean and
variance are representable as linear functionals of a common
unobserved Poisson random measure. The analysis of these models is
facilitated by applying the methods in James~(2005, 2002), in
particular an Esscher type transform of Poisson random measures;
in conjunction with a special case of the Weber-Sonine formula. It
is shown that the marginal likelihood may be expressed in terms of
a multidimensional Fourier-cosine transform. This yields tractable
forms of the likelihood and also allows a full Bayesian posterior
analysis of the integrated volatility process. A general formula
for the posterior density of the log price given the observed data
is derived, which could potentially have applications to option
pricing. We also identify tractable subclasses, where inference
can be based on a finite number of independent random variables.
We close by obtaining explicit expressions for likelihoods
incorporating leverage. It is shown that inference does not
necessarily require simulation of random measures. Rather,
classical numerical integration can be used in the most general
cases.
 \EndQuote
%\baselineskip14pt
%\begin{document}
\rm
%\newpage

\section{Introduction}
Barndorff-Nielsen and Shephard~(2001a, b)(BNS) introduce a class
of continuous time stochastic volatility~(SV) models that allows
for more flexibility over Gaussian based models such as the
Black-Scholes model[see Black and Scholes~(1973) and
Merton~(1973)]. Their proposed SV model is based on the following
differential equation, \Eq dx^{*}(t)=(\mu+\beta
v(t))dt+v^{1/2}(t)dw(t) \label{BNS}\EndEq where $x^*(t)$ denotes
the log-price level, $w(t)$ is Brownian motion, and independent of
$w(t)$, $v(t)$ is a stationary Ornstein-Uhlenbeck (OU) process
which models the {\it instantaneous volatility}. This model is an
extension of the Black-Scholes or Samuelson model which arises by
replacing $v$ with a fixed variance, say $\sigma^{2}$.  The
additional innovation in BNS is that modeling volatility as a
random process, $v(t)$, rather than a random variable, not only
allows for heavy-tailed models, but additionally induces serial
dependence. This serial dependence is used to account for a
clustering affect referred to as {\it volatility persistence}. The
work of Carr, Geman, Madan, and Yor~(2003) discuss this point
further. See also Duan~(1995) and Engle~(1982) for different
approaches to this type of phenomenon. The model of BNS has gained
a great deal of interest with some related works including Carr,
Geman, Madan, and Yor~(2003), Barndorff-Nielsen and
Shephard~(2003), Eberlein~(2001), Nicolato and Venardos~(2001),
Benth, Karlsen, and Reikvam~(2003). See also the discussion
section in~Barndorff-Nielsen and Shephard~(2001a).

One current drawback of this approach is the unavailability of a
tractable analysis of likelihood based stochastic volatility
models for the returns of log prices of stocks. These models are
based on the integrated volatility process
$\tau(t)=\int_{0}^{t}v(u)du$. Several MCMC procedures have been
proposed to handle subclasses of these models requiring simulation
of points from random processes. See for instance, Roberts,
Papaspiliopoulos and Dellaportas~(2004) and the discussion section
in Barndorff-Nielsen and Shephard~(2001a).

In this paper, we shall actually provide a complete analysis of a
significantly more complex class of likelihood models.
Specifically models where $\tau$ is expressible as a linear
functional of a Poisson random measure. This includes the
superposition processes mentioned in Barndorff-Nielsen and
Shephard~(2001a, b) and much more general spatial models for
$\tau$. Our results are therefore applicable to a wide range of
models and applications. We first present a description of these
models similar to the framework outlined in Barndorff-Nielsen and
Shephard~(2001a). We then show how the results in James~(2005,
2002) are easily applied to this setting via the usage of a Bessel
integral identity involving the cosine function. That is, a
special case of what is called the Weber-Sonine formula. This
leads to an interesting series of tractable characterizations of
such processes, including an identification of simple subclasses
of these models. As a byproduct, we derive the posterior
predictive density of what could be considered as models for the
log price of stocks. This may prove useful to applications in
option pricing. Moreover, our methods do not require simulation of
random measures and in the most general cases can be handled by
more classical numerical integration methods. We point out also
that procedures to simulate from random measures often require
explicit knowledge of the L\'evy density associated with an
infinitely divisible random variable. This is an important point
as there are some interesting cases where the probability density
of a random variable is known explicitly but its corresponding
L\'evy density is unknown. Our results show that one only requires
knowledge of the form of the L\'evy exponent or log of the Laplace
transform of a corresponding random variable. Inference using the
classes of models that we identify in sections 2.5 and 4 can be
performed based on at most $J<\infty$ independent {\it latent}
random variables. In the case of section 4, $J=n+1$, where $n$
denotes the number of observations. Section 2.6 discusses another
class where exact calculations of a different form are easily
obtained. Section 5 describes exact expression for likelihoods of
generalized types of leverage effects models.

\Remark The appearance of integrals involving Bessel functions is
certainly not new to applications in finance as can be seen in the
case of the important work of Yor~(1992) on Asian Options. See
also Carr and Schr\"oder~(2004). We shall however employ a
different, but certainly related, integral identity. \EndRemark
\subsection{Likelihood model and representation}
The model of Barndorff-Nielsen and Shephard~(2001a, section 5.4)
translates into a likelihood based model as follows.
 Let ${X_{i}}$ for $i=1,\ldots, n$ denote a sequence of aggregate
returns of the log price of a stock observed over intervals of
length $\Delta>0$. Additionally for each interval
$[(i-1)\Delta,i\Delta]$, let
$\tau_{i}=\tau(i\Delta)-\tau((i-1)\Delta)$. Now the model in
~\mref{BNS} implies that $X_{i}|\tau_{i},\beta, \mu$ are
conditionally independent with \Eq
X_{i}=\mu\Delta+\tau_{i}\beta+\tau^{1/2}_{i}\epsilon_{i},
\label{data}\EndEq where $\epsilon_{i}$ are independent standard
Normal random variables. Hence if $\tau$ depends on external
parameters $\theta$, one is interested in estimating
$(\mu,\beta,\theta)$ based on the likelihood \Eq
\Lcr(\X|\mu,\beta,\theta)=\int_{{\mathbb R}^{n}_{+}}
\[\prod_{i=1}^{n}\phi(X_{i}|\mu\Delta+\beta
\tau_{i},\tau_{i})\]f(\tau_{1},\ldots,\tau_{n}|\theta)d\tau_{1},\ldots,d\tau_{n}
\label{BNSlik} \EndEq where, setting $A_{i}=(X_{i}-\mu\Delta)$,
and ${\bar A}=n^{-1}\sum_{i=1}^{n}A_{i}$,
$$\phi(X_{i}|\mu\Delta+\beta
\tau_{i},\tau_{i})={\mbox e}^{A_{i}\beta
}\frac{1}{\sqrt{2\pi}}\tau^{-1/2}_{i}{\mbox
e}^{-A^{2}_{i}/(2\tau_{i})}{\mbox e}^{-\tau_{i}\beta^{2}/2}
$$
denotes a Normal density.

The quantity $f(\tau_{1},\ldots, \tau_{n}|\theta)$ denotes the
joint density of the integrated volatility based on the intervals
$[(i-1)\Delta,i\Delta]$ for $i=1,\ldots, n$. Barndorff-Nielsen and
Shephard~(2001a) note that the likelihood is intractable and hence
makes exact inference difficult. The apparent intractability is
attributed to the complex nature of $f(\tau_{1},\ldots,
\tau_{n}|\theta)$ which is derived from a random measure. However,
we shall show that in fact it is quite easy to deal with
$f(\tau_{1},\ldots, \tau_{n}|\theta)$ for more general $\tau$ by
means of the Poisson partition calculus methods outlined in
James~(2005, 2002). Rather, the stumbling block which currently
prevents one from integrating out the infinite-dimensional
components in the likelihood, is inherent from the Normal
distribution of $X_{i}|\tau_{i},\beta,\mu$. Quite simply the
Normal assumption yields exponential terms of the form
$$
{\mbox e}^{-A^{2}_{i}/(2\tau_{i})}{\mbox { rather than }}{\mbox
e}^{-\tau_{i}A^{2}_{i}}.
$$
In the next section we apply an integral identity to circumvent
this problem.
\subsection{Bessel integral representation of the likelihood}
 In order to calculate ~\mref{BNSlik} we first employ a
Bessel integral identity which we state in more general terms.
Suppose that $J_{v}(x)$ denotes a Bessel function of the first
kind of order $v$. Then for $v>-1$, and numbers $a, p$ \Eq
{p}^{-2(v+1)}{\mbox
e}^{-a^2/4p^{2}}=2^{v+1}a^{-v}\int_{0}^{\infty}J_{v}(at)t^{v+1}{\mbox
e}^{-p^{2}t^{2}}dt. \label{key1}\EndEq This is a special case of
the Weber-Sonine formula. See for instance Andrews, Askey and
Roy~(1999, p.222) and Watson~(1966, p. 394 eq. (4)) for the
identity and also those references for Bessel functions. Taking
for each $i$, $p^{2}=\tau_{i}/2$, $a=|A_{i}|$ and $v=-1/2$, it
follows the marginal likelihood is given by, \Eq
\Lcr(\X|\mu,\beta,\theta)=\frac{{\mbox e}^{n{\bar A}\beta
}}{\pi^{n}}\int_{{\mathbb R}^{n}_{+}} \E\[\prod_{i=1}^{n}{\mbox
e}^{-(y_{i}^{2}/2+\beta^{2}/2)\tau_{i}}\]
\prod_{i=1}^{n}\cos(y_{i}|A_{i}|)dy_{i}\label{key2} \EndEq where
$$\E\[\prod_{i=1}^{n}{\mbox
e}^{-(y_{i}^{2}/2+\beta^{2}/2)\tau_{i}}\]=\int_{{\mathbb
R}^{n}_{+}} \prod_{i=1}^{n}{\mbox
e}^{-(y_{i}^{2}/2+\beta^{2}/2)\tau_{i}}f(\tau_{1},\ldots,\tau_{n}|\theta)\prod_{l=1}^{n}d\tau_{l}.$$
\Remark Notice that the expression~\mref{key2} has nothing to do
with the distributional properties of $\tau$. The appearance of
the cosine in~\mref{key2} is due to the identity $$
\sum_{k=0}^{\infty}
\frac{{(-1)}^{k}x^{2k}}{\Gamma(k+1/2)k!2^{2k-1/2}}=x^{1/2}J_{-1/2}(x)=\sqrt{\frac{2}{\pi}}\cos(x)
$$ where $J_{-1/2}(x)$ is a Bessel function of the first kind of
order $-1/2$. See for example Andrews, Askey and Roy~(1999, p.
202). In other words we are using the classical Fourier-Cosine
Integral
$$
\frac{1}{\pi}\int_{0}^{\infty}\cos(y|A_{i}|){\mbox
e}^{-\frac{y^{2}\tau_{i}}{2}}dy=\frac{1}{\sqrt{2\pi}}\tau^{-1/2}_{i}{\mbox
e}^{-\frac{A^{2}_{i}}{2\tau_{i}}}
$$
\EndRemark

\section{Evaluation of the likelihood for general $\tau$}
The representations in ~\mref{key2} and~\mref{key1} allows us to
immediately apply the results in James~(2005) to obtain a full
analysis for quite general $\tau$, which we now describe. Let $N$
denote a Poisson random measure on some Polish space $\Vcr$ with
mean intensity,
$$
\E[N(dx)|\nu]=\nu(dx).
$$
We denote the Poisson law of $N$ with intensity $\nu$ as
$\Pe(dN|\nu)$. The Laplace functional for $N$ is defined as
$$
\E[{\mbox e}^{-N(f)}|\nu]=\int_{\M}{\mbox
e}^{-N(f)}\Pe(dN|\nu)={\mbox e}^{-\Lambda(f)}
$$
where for any positive $f$, $N(f)=\int_{\Vcr}f(x)N(dx)$ and
$\Lambda(f)=\int_\Vcr(1-{\mbox e}^{-f(x)})\nu(dx).$ $\M$ denotes
the space of boundedly finite measures on $\Vcr$ [see Daley and
Vere-Jones~(1988)].
 We suppose that
$\tau_{i}=N(f_{i})$, for $i=1,\ldots, n$ where
$f_{1},\ldots,f_{n}$ are positive measureable functions on $\Vcr$.
Notice now that the index $i=1,\ldots,n$ need not correspond to
fixed intervals involving $\Delta$. With this in mind, let
$(w_{1},\ldots,w_{n})$ denote arbitrary non-negative numbers.
Define for $i=1,\ldots, n$, functions
$R_{i}(x)=\sum_{j=1}^{i}w_{j}f_{j}(x)$ and $\nu_{R_{i}}(dx)={\mbox
e}^{-R_{i}(x)}\nu(dx).$ Then all our results will follow from the
following special case of James~(2005, Proposition 2.1), which can
be viewed as an {\it Esscher}-type transform,\Eq {\mbox
e}^{-N(\sum_{i=1}^{n}w_{i}f_{i})}\Pe(dN|\nu)=\Pe(dN|\nu_{R_{n}}){\mbox
e}^{-\Lambda(\sum_{i=1}^{n}w_{i}f_{i})}.\label{esscher}\EndEq
Additionally the following decomposition is sometimes useful
$$ \E\[{\mbox e}^{-N(\sum_{i=1}^{n}w_{i}f_{i})}|\nu\] ={\mbox
e}^{-\Lambda(\sum_{i=1}^{n}w_{i}f_{i})}=\E\[{\mbox
e}^{-N(w_{1}f_{1})}|\nu\]\prod_{i=2}^{n} \E\[{\mbox
e}^{-N(w_{i}f_{i})}|\nu_{R_{i-1}}\].$$ This expression appears in
James~(2002) and may be  obtained by repeated application
of~\mref{esscher}.

Now, throughout, for each $n\ge 1$, define $
\Omega_{n}(x)=\sum_{i=1}^{n}(y^{2}_{i}/2+\beta^{2}/2)f_{i}(x). $
This is a special case of $\sum_{i=1}^{n}w_{i}f_{i}$ with
$w_{i}=(y^{2}_{i}/2+\beta^{2}/2)$ for $i=1,\ldots,n$. The
following result is immediate from an application of Fubini's
theorem,~\mref{esscher} and the representation~\mref{key2}.
\begin{thm} Suppose that $\tau_{i}=N(f_{i})$ for $i=1,\ldots,n$
where $N$ is a Poisson random measure on $\Vcr$ with intensity
$\nu$. Then setting $w_{i}=y_{i}^{2}/2+\beta^{2}/2$ for
$i=1,\ldots,n,$ in ~\mref{esscher}, the likelihood~\mref{BNSlik}
can be expressed as $$\Lcr(\X|\mu,\beta,\theta)= \frac{{\mbox
e}^{n{\bar A}\beta }}{\pi^{n}}\int_{{\mathbb R}^{n}_{+}} {\mbox
e}^{-\Lambda(\Omega_{n})}
\prod_{i=1}^{n}\cos(y_{i}|A_{i}|)dy_{i}.$$ \qed\end{thm}

\subsection{Posterior distribution of parameters}
Theorem 2.1 shows that Bayesian inference for $(\mu,\beta,\theta)$
may be described as follows.
\begin{prop} Suppose that $\tau$ depends on a $d$-dimensional parameter
$\theta$. Then if $q(d\theta)$, $q(d\beta)$, $q(d\mu)$ denote
independent prior distributions for $(\beta,\mu,\theta)$, their
posterior distribution can be written as, $$
q(d\beta,d\mu,d\theta|\X)\propto \int_{{\mathbb R}^{n}_{+}}
[{\mbox e}^{-[\Lambda_{\theta}(\Omega_{n})-n{\bar A}\beta
]}q(d\theta)q(d\beta)]
\prod_{i=1}^{n}\cos(y_{i}|A_{i}|)dy_{i}q(d\mu),$$ where
$\Lambda_{\theta}$ denotes the dependence of $\Lambda$ on
$\theta.$\qed
\end{prop}

\subsection{Posterior distribution of the process}
The above results describe the behaviour of the finite-dimensional
likelihood and parameters. It is useful to also obtain a
description of the underlying random process given the data. This
allows one to see directly how the data affects the overall
process. Moreover, combined with the results in James~(2005), it
provides a calculus for more general functionals. Define the
measure,
$$\Qcr_{\X}(d\y)=\pi^{-n}{\mbox e}^{-[\Lambda(\Omega_{n})-n{\bar
A}\beta]} \prod_{i=1}^{n}\cos(y_{i}|A_{i}|)dy_{i}
/\Lcr(\X|\mu,\beta,\theta).
$$
For notational simplicity we suppose that $(\mu,\beta,\theta)$ are
fixed. The next result also follows immediately from an
application of Fubini's theorem,~\mref{esscher} and the
representation~\mref{key2}.
\begin{thm} Suppose that the distribution of  $\X$  is given by
~\mref{BNSlik}, and that $\tau$ and $N$ are defined by the
specifications in Theorem 2.1. Let $\Omega_{n}(x)=\sum_{i=1}
(y_{i}^{2}/2+\beta^{2}/2)f_{i}(x)$. Then the posterior
distribution of $N|\X$ is given by the mixture $$ \int_{{\mathbb
R}^{n}_{+}}\Pe(dN|\nu_{\Omega_{n}})\Qcr_{\X}(d\y)$$ which
determines the posterior distribution of $\tau$ and related
quantities. $\Pe(dN|\nu_{\Omega_{n}})$ can be viewed as the
posterior distribution of $N$ given the information in $\Y,\X$ and
corresponds to the law of a Poisson random measures with mean
intensity $$\nu_{\Omega_{n}}(dx):={\mbox
e}^{-\Omega_{n}(x)}\nu(dx)=\nu(dx){\mbox e}^{-\sum_{i=1}
(y_{i}^{2}/2+\beta^{2}/2)f_{i}(x)}.$$
 \qed
\end{thm}
The next result, which gives an expression for the posterior
Laplace functional of $N$, is an immediate consequence
of~\mref{esscher} combined with Theorem 2.2.
\begin{prop} The posterior Laplace functional of $N|\X$, according
to the Theorem 2.2, is given by
$$
\int_{{\mathbb R}^{n}_{+}}\[\int_{\M}{\mbox
e}^{-N(f)}\Pe(dN|\nu_{\Omega_{n}})\]\Qcr_{\X}(d\y)=\int_{{\mathbb
R}^{n}_{+}}{\mbox
 e}^{-[\Lambda(f+\Omega_{n})-\Lambda(\Omega_{n})]}\Qcr_{\X}(d\y)
 $$
for $f$ such that $\Lambda(f+\Omega_{n})<\infty.$ \qed\end{prop}
\subsection{A general posterior predictive density for the log price}
We now define a random variable similar to~\mref{data} which can
be thought of as representing the log-price and give an explicit
expression for its posterior density given $\X$. The random
variable is defined as, \Eq {\tilde X}=\mu{\tilde \Delta}+{\tilde
\tau}\beta+{{\tilde \tau}}^{1/2}{\tilde \epsilon}
\label{abprice}\EndEq where ${\tilde \Delta}$ denotes a general
positive quantity, ${\tilde \epsilon}$ is a standard Normal random
variable independent of all other variables, and ${\tilde
\tau}=N({\tilde f})$ for some positive function ${\tilde f}$ such
that Laplace transform of ${\tilde \tau}$ exists. \begin{prop}
Suppose $N$ and the data structure of $\X$ is defined as in
Theorem 2.2. Let ${\tilde X}$ be defined by~\mref{abprice}. Denote
its marginal density as $f_{{\tilde X}}(\cdot|\beta,\mu)$ and its
posterior density given the data $\X$ from~\mref{data} as
$f_{{\tilde X}}(\cdot|\beta,\mu,\X)$. Then the following results
hold\Enumerate
\item[(i)]$f_{{\tilde X}}(x|\beta,\mu)=
\frac{1}{\pi}{{\mbox e}^{{(x-{\mu\tilde \Delta})}\beta}}
\int_{0}^{\infty} {\mbox e}^{-\Lambda((y^{2}/2+\beta^{2}/2){\tilde
f})}\cos(y|x-\mu{\tilde \Delta}|)dy.$
\item[(ii)]The posterior density of the log stock price given $\X$
is, $f_{{\tilde X}}(x|\beta,\mu,\X)$, given by
$$\frac{1}{\pi}{{\mbox e}^{{(x-{\mu\tilde \Delta})}\beta}}
\int_{{\mathbb R}^{n}_{+}}\[\int_{0}^{\infty}{\mbox
e}^{-[\Lambda(w_{n+1}{\tilde
f}+\Omega_{n})-\Lambda(\Omega_{n})]}\cos(y|x-\mu{\tilde
\Delta}|)dy\]\Qcr_{\X}(d\y),$$ where
$w_{n+1}=(y^{2}/2+\beta^{2}/2).$\qed \EndEnumerate
\end{prop}
\Proof Setting ${\tilde
\Omega}_{n+1}(x)=\Omega_{n}(x)+w_{n+1}{\tilde f}(x)$, the results
follow from~\mref{key1} and Theorem 2.2, using the fact
from~\mref{esscher} that,
$${\mbox e}^{-N(w_{n+1}{\tilde f})}\Pe(dN|\nu_{\Omega_{n}})=\Pe(dN|\nu_{{\tilde
\Omega}_{n+1}}){\mbox e}^{-[\Lambda(w_{n+1}{\tilde
f}+\Omega_{n})-\Lambda(\Omega_{n})]}.
$$
\EndProof

\subsection{Simplifications for a class of $\tau$ via an inversion formula}
We have shown that for $\tau$ modeled quite generally that its
contribution to the likelihood~\mref{BNSlik} is only through the
exponent $\Lambda(\sum_{i=1}^{n}w_{i}f_{i})$. That is through the
form of $\sum_{i=1}^{n}w_{i}f_{i}$ and $\nu$. With a view towards
choosing $\tau$ which are the most tractable we present the
following interesting result.

\begin{thm} Suppose that for arbitrary non-negative
$(w_{1},\ldots,w_{n})$, and an integer $J$, there is an array of
non-negative numbers $(a_{ij})$ such that
$\Lambda(\sum_{i=1}^{n}w_{i}f_{i})=\sum_{j=1}^{J}\Lambda([\sum_{i=1}^{n}w_{i}a_{ij}]h_{j}).$
Where $(h_{j})$ are non-negative functions on $\Vcr$ such that
$\Lambda(\omega h_{j})<\infty$ for all $\omega\ge 0$. Let
$(T_{1},\ldots, T_{J})$ denote $J$ independent random variables
with respective Laplace transforms $E[{\mbox e}^{-\omega
T_{j}}]={\mbox e}^{-\Lambda(\omega h_{j})}$ for $j=1,\ldots,J$.
Moreover $\Lcr(\X|\mu,\beta,\theta)$ denotes the likelihood given
in~\mref{BNSlik}. Then, \Enumerate
\item[(i)]$\E\[{\mbox e}^{-\sum_{i=1}^{n}w_{i}\tau_{i}}\]=\prod_{j=1}^{J}{\mbox
e}^{-\Lambda([\sum_{i=1}^{n}w_{i}a_{ij}]h_{j})}=
\prod_{j=1}^{J}E[{\mbox e}^{-[\sum_{i=1}^{n}w_{i}a_{ij}]T_{j}}].$
\item[(ii)]$\Lcr(\X|\mu,\beta,\theta)=
\E\[\prod_{i=1}^{n}\phi(X_{i}|\mu\Delta+\beta
[\sum_{j=1}^{J}a_{ij}T_{j}],[\sum_{j=1}^{J}a_{ij}T_{j}])\],$ where
the expectation is respect to the distribution of
$(T_{1},\ldots,T_{J}).$\EndEnumerate \qed
\end{thm}
\Proof Statement (i) is immediate from the specification of
$\Lambda(\sum_{i=1}^{n}w_{i}f_{i})$. Statement (i) implies that
that one may replace $\sum_{i=1}^{n}w_{i}\tau_{i}$ with
$\sum_{i=1}^{n}w_{i}[\sum_{j=1}^{J}a_{ij}T_{j}]$. Now setting each
$w_{i}= (y_{i}^{2}/2+\beta^{2}/2)$ for $i=1,\ldots,n$, one uses
~\mref{key1} and~\mref{key2} to conclude the result. \EndProof
Statement (ii) of Theorem 2.3 allows one to approximate the
likelihood by the simulation of $J$ independent random variables.
It also demonstrates that it is rather straightforward to conduct
{\it parametric} Bayesian or frequentist estimation procedures,
where $(T_{1},\ldots,T_{J})$ are viewed as independent latent
variables. The next proposition puts this in a Bayesian framework.
\begin{prop} Suppose that $(T_{1},\ldots, T_{J})$ depend on
external parameters, say $\theta$. Then assuming a joint prior
$q(d\theta,d\mu,d\beta)$, posterior inference may be obtained
based on the model derived from augmenting the likelihood in
Theorem 2.3. That is, the joint distribution of
$(\X,T_{1},\ldots,T_{J},\theta,\mu,\beta)$ given by
$$
\[\prod_{i=1}^{n}\phi(X_{i}|\mu\Delta+\beta
[\sum_{j=1}^{J}a_{ij}T_{j}],[\sum_{j=1}^{J}a_{ij}T_{j}])\]
\[\prod_{j=1}^{J}f_{T_{j}}(T_{j})\]q(d\theta,d\mu,d\beta).
$$
\qed
\end{prop}
\section{Tractable expressions}
It is noted that, at first glance, one may find it difficult to
work with the expressions involving cosines. Here, influenced by
some arguments in Devroye~(1986a), we give a representation of the
likelihood that can be numerically evaluated via the simulation of
random variables. First let ${\bf p}=\{p_{1},\ldots,p_{n}\}$
denote a vector of positive numbers and for each $i$,let
$$
H(y_{i}|p_{i})=\frac{2}{\sqrt{2\pi p_{i}}}{\mbox
e}^{-\frac{y^{2}_{i}}{2p_{i}}}{\mbox { for }}y_{i}>0
$$
denote a half normal density.  Now, notice that $ 0\leq
1-\prod_{i=1}^{n}\cos(y_{i})\leq 2$, and \Eq \int_{{\mathbb
R}^{n}_{+}}\[1-\prod_{i=1}^{n}\cos(y_{i}|A_{i}|)\]
H(y_{i}|p_{i})dy_{i}=1-{\mbox e}^{-\frac{\sum_{i=1}^{n}A^{2}_i
p_{i}}{2}}=C_{n}({\mathbf A},{\mathbf p})\label{cosden} \EndEq
From these facts we describe a joint density
\begin{prop}Augmenting the expression in~\mref{cosden} leads to a
joint density of an array of positive random variables
$\Y=\{Y_{1,n},\ldots,Y_{n,n}\}$ given by,
$$
r_{n}(\y|{\bf p})=\frac{\[1-\prod_{i=1}^{n}\cos(y_{i}|A
_{i}|)\]\prod_{i=1}^{n}H(y_{i}|p_{i})}{ C_{n}({\mathbf A},{\mathbf
p})}
$$
Equivalently, for $k=1,\ldots, n$, the conditional density of
$Y_{k,n}|Y_{1,n},\ldots, Y_{k-1,n}$ is proportional to
$[1-\lambda_{k}cos(y_{k}|A_{k}|)]H(y_{k}|p_{k})$, where
$\lambda_{k}={\mbox
e}^{-\sum_{i=k+1}^{n}\frac{A^{2}_{i}p_{i}}{2}}\prod_{i=1}^{k-1}\cos(y_{i}|A_{i}|)$
for $k=2,\ldots, n-1$, $\lambda_{1}= {\mbox
e}^{-\sum_{i=2}^{n}\frac{A^{2}_{i}p_{i}}{2}}$, and
$\lambda_{n}=\prod_{i=1}^{n-1}\cos(y_{i}|A_{i}|).$
\end{prop}
Define the function
$$
\Upsilon_{n}(\beta,\theta):=\frac{1}{\pi^{n}}\int_{{\mathbb
R}^{n}_{+}} {\mbox e}^{-\Lambda(\Omega_{n})}
\prod_{i=1}^{n}dy_{i}=\E\[\prod_{i=1}^{n}\frac{{\mbox
e}^{-\beta^{2}\tau_{i}}}{\sqrt{2\pi \tau_{i}}}\]\leq
\E\[\prod_{i=1}^{n}\frac{1}{\sqrt{\tau_{i}}}\]
$$

These points lead to following representation of the likelihood.
\begin{prop} Suppose that for fixed $n$,
$\E\[\prod_{i=1}^{n}\frac{1}{\sqrt{\tau_{i}}}\]<\infty$, then the
likelihood in Theorem 2.1 may be written as
$$
{\mbox e}^{{\bar A}\beta
}\[\Upsilon_{n}(\beta,\theta)-\frac{C_{n}({\mathbf A},{\mathbf
p})}{\pi^{n}}\E\[\frac{{\mbox
e}^{-\Lambda(\Omega_{n})}}{\prod_{i=1}^{n}H(Y_{i,n}|p_{i})}\]\]
$$
where the random vector $\{Y_{1,n},\ldots,Y_{n,n}\}$ has its joint
distribution described by proposition 3.1, and
$\Omega_{n}(x)=\sum_{i=1}^{n}[(Y^{2}_{i,n}+\beta^{2})/2]g_{i}(x)$
\end{prop}
\Remark Proposition 3.2 shows that one may approximate the
likelihood by simulating random variables described in Proposition
3.1. Such an approach should work well with a Bayesian procedure
for estimating the parameters $(\mu,\beta,\theta)$. Methods to
easily sample the random variables in proposition 3.1, may be
deduced from Devroye~(1986a, b). Alternatively one may sample from
the densities $H(y_{i}|p_{i})$. One may also use other
densities.\EndRemark

\section{Analysis of the BNS-OU model}
In this section we will show how our results apply to the basic
integrated volatility model of Barndorff-Nielsen and
Shephard~(2001a, b). We shall refer to this model as the BNS-OU
model. First suppose that $N$ is a Poisson random measure on
$(0,\infty)\times(-\infty,\infty)$ with intensity
$$
\nu(du,dy)=\rho(du)dy
$$
where $\rho$ is the L\'evy density of an infinite-divisible random
variable, say $T$, with Laplace transform for $\omega\ge 0$,
$$
\E[{\mbox e}^{-\omega T}]={\mbox e}^{-\psi(\omega)} {\mbox { where
}}
 \psi(\omega)=\int_{0}^{\infty}(1-{\mbox e}^{-\omega
u})\rho(du).
$$
Now we model the background driving L\'evy process (BDLP), say
$z$, as a completely random measure which is expressible in
distribution as $ z(dt)=\int_{0}^{\infty}uN(du,dt). $ Note that
for any non-negative function $g$ on $(-\infty,\infty)$, it
follows that
$$
z(g)=\int_{-\infty}^{\infty}g(y)z(dy)=N(f_{g})$$ where
$f_{g}(u,y)=ug(y)$ on $(0,\infty)\times(-\infty,\infty)$.
 Additionally,
$$\E\[{\mbox e}^{-z(g)}|\nu\]=\int_{\M}{\mbox
e}^{-N(f_{g})}\Pe(dN|\nu)={\mbox
e}^{-\int_{-\infty}^{\infty}\psi(g(y))dy}={\mbox
e}^{-\Lambda(f_{g})}.$$  One may express the Barndorff-Nielsen and
Shephard~(2001 a, b) integrated OU process $\tau$ as \Eq
\tau(t)=\lambda^{-1}[(1-{\mbox e}^{-\lambda
t})\int_{-\infty}^{0}{\mbox e}^{y}z(dy)+\int_{0}^{t}(1-{\mbox
e}^{-\lambda(t-y)})z(dy)]\label{model1} \EndEq where
$v(0):=v_{0}=\int_{-\infty}^{0}{\mbox e}^{y}z(dy)$. The form in
~\mref{model1} is taken from Carr, Geman, Madan and Yor~(2003, p.
365). It follows that for any $s<t$,
$[\tau(t)-\tau(s)]=z(g_{s,t})=N(f_{s,t})$ where
$f_{s,t}(u,y)=ug_{s,t}(y)$ and $\lambda g_{s,t}(y)$ equals, \Eq
{\mbox e}^{-\lambda s}(1-{\mbox e}^{-\lambda (t-s)}){\mbox
e}^{y}I_{\{y\leq 0\}}+(1-{\mbox e}^{-\lambda(t-y)})I_{\{s<y \leq
t\}}+{\mbox e}^{-\lambda s}(1-{\mbox e}^{-\lambda(t-s)}){\mbox
e}^{\lambda y}I_{\{0<y\leq s\}}.\label{gendif3}\EndEq The first
component in~\mref{gendif3} represents the contribution from
$v_{0}$. Specializing this to $s=(i-1)\Delta$ and $t=i\Delta$ one
has $\tau_{i}=z(g_{i,1}+g_{i,2})=N(f_{i})$ where
$f_{i}(u,y)=u[g_{i,1}(y)+g_{i,2}(y)]$ and \Eq
g_{i,1}(y)=\lambda^{-1}[(1-{\mbox e}^{-\lambda
(i\Delta-y)})I_{\{(i-1)\Delta<y\leq i\Delta\}}+{\mbox
e}^{-\lambda(i-1)\Delta}(1-{\mbox e}^{-\lambda\Delta}){\mbox
e}^{y}I_{\{y\leq 0\}}]\label{gform} \EndEq and \Eq
g_{i,2}(y)=\lambda^{-1}{\mbox e}^{-\lambda(i-1)\Delta}(1-{\mbox
e}^{-\lambda \Delta}){\mbox e}^{\lambda y}I_{\{0<y\leq
(i-1)\Delta\}}.\label{hform}\EndEq Now for $i=1,\ldots, n$, set
$r_{i}=\lambda^{-1}[\sum_{k=i}^{n}w_{k}{\mbox e}^{-\lambda (k-1)
\Delta}](1-{\mbox e}^{- \lambda \Delta})$. Now notice that for any
sequence of numbers, the simplest expression will be obtained by
utilizing the following facts. \Eq
\sum_{j=1}^{n}w_{j}[g_{j,1}{(y)}+g_{j,2}(y)]=r_{1}{\mbox
e}^{y}{\mbox { for }}y\leq 0 \label{id1}\EndEq and for
$i=1,\ldots, n$\Eq
\sum_{j=1}^{n}w_{j}[g_{j,1}{(y)}+g_{j,2}(y)]=\zeta
(y|w_{i},r_{i+1}){\mbox { for }}(i-1)\Delta<y\leq i\Delta.
\label{id2}\EndEq Where for each $i$, $\zeta
(y|w_{i},r_{i+1})=[\lambda^{-1}w_{i}(1-{\mbox
e}^{-\lambda(i\Delta-y)})+r_{i+1}{\mbox e}^{\lambda y}].$

\begin{prop} For $0\leq s<t$, let $\tau(t)-\tau(s)$ be defined by~\mref{model1} and \mref{gendif3}
Then the results of Theorem 2.1 and 2.2 hold with
$f_{i}(u,y)=u[g_{i,1}(y)+g_{i,2}(y)]$,
 $w_{i}=(y^{2}_{i}/2+\beta^{2}/2)$, as described in~\mref{gform}
and~\mref{hform}. In particular, using a change of variable,
\Enumerate
\item[(i)]
$ {\mbox e}^{-\Lambda(\sum_{i=1}^{n}w_{i}f_{i})}={\mbox
e}^{-\Phi_{0}(r_{1})}{\mbox
e}^{-\Phi_{n}(w_{n})}\prod_{i=1}^{n-1}{\mbox
e}^{-\Phi_{i}(w_i|r_{i+1})} $ \item[(ii)]
$\Phi(w_{i}|r_{i+1})=\int_{1-{\mbox
e}^{-\lambda\Delta}}^{1}\lambda^{-1}\psi(r_{i+1}{\mbox e}^{\lambda
i\Delta}(1-u)+\lambda^{-1}w_{i}u)\frac{du}{1-u},$ for $i=1,\ldots,
n-1$
\item[(iii)]
$\Phi(w_{n})=\int_{1-{\mbox
e}^{-\lambda\Delta}}^{1}\lambda^{-1}\psi(\lambda^{-1}w_{n}u)\frac{du}{1-u}$
\item[(iv)]
$\Phi_{0}(r_{1})=\int_{0}^{1}\psi(r_{1}u)\frac{du}{u}$, where
${\mbox e}^{-\Phi_{0}(r_{1})}=\E[{\mbox e}^{-r_{1}v_{0}}].$
\EndEnumerate \qed
\end{prop}
\Remark Expressions of the form in [(iii)] of Proposition 3.1 are
known to be a key component in option pricing using the BNS-OU
model. However explicit calculations have only been given for a
few cases. See Barndorff-Nielsen and Shephard~(2003), Nicolato and
Venardos~(2003) and Carr, Geman, Madan and Yor~(2003). Note that
if for $y>0$, we change the Lebesque measure,
 $dy$, to ${\mbox e}^{\lambda y}dy$, the calculations for
$\Phi(w_{i}|r_{i+1})$ for $i=1,\ldots, n,$ where
$\Phi(w_{n}|r_{n+1})=\Phi(w_{n})$, are greatly simplified. See
section 4.2 for a closely related discussion. \EndRemark

\section{Analysis of a simple class of models}
This last section, which is based on a class of models from
section 2.4, examines models which are the most tractable and we
believe still flexible enough to be applied to general classes of
problems. Implicitly, we are taking a Bayesian nonparametrics
viewpoint of seeking random measures as priors which are both
flexible in a modeling sense and easily manipulated. For
concreteness, we start out with a variation of the
Barndorff-Nielsen and Shephard~(2001a,b, 2003) integrated OU
process $\tau$. Here we set, \Eq\tau(t)=\lambda^{-1}[(1-{\mbox
e}^{-\lambda t})\int_{-\infty}^{0}{\mbox
e}^{y}z(dy)+\int_{0}^{\lambda t}(1-{\mbox
e}^{-y})z(dy)]\label{model2},\EndEq where again
$v(0)=v_{0}=\int_{-\infty}^{0}{\mbox e}^{y}z(dy)$. Interestingly
from Barndorff-Nielsen and Shephard~(2003, p. 282), one has the
following distributional equivalence of marginal distributions,
$$
\lambda^{-1}\int_{0}^{\lambda t}(1-{\mbox e}^{-y})z(dy)\overset
{d}=\lambda^{-1}\int_{0}^{\lambda t}(1-{\mbox
e}^{-\lambda(t-y)})z(d\lambda y)
$$
where the right hand side equates with the model in section 3.
However, now for $s<t$, $$ \lambda[\tau(t)-\tau(s)]={\mbox
e}^{-\lambda s}(1-{\mbox e}^{-\lambda
(t-s)})\int_{-\infty}^{0}{\mbox e}^{y}z(dy)+\int_{0}^{\lambda
t}(1-{\mbox e}^{-y})I_{\{\lambda s<y \leq \lambda t\}}z(dy).$$
Specializing this to $t=i\Delta$ and $s=(i-1)\Delta$ yields
$$
\tau_{i}=\int_{-\infty}^{\infty}g_{i}(y)z(dy)=z(g_{i})=N(f_{i})$$
where $f_{i}(u,y)=ug_{i}(y)$ with
$$g_{i}(y)=\lambda^{-1}[(1-{\mbox
e}^{-y})I_{\{\lambda(i-1)\Delta<y\leq \lambda i\Delta\}}+{\mbox
e}^{-\lambda(i-1)\Delta}(1-{\mbox e}^{-\lambda\Delta}){\mbox
e}^{y}I_{\{y\leq 0\}}].$$ Now notice that for any sequence of
numbers, the simplest expression will be obtained by utilizing the
following facts. \Eq
\sum_{j=1}^{n}w_{j}g_{j}{(y)}=\lambda^{-1}[\sum_{i=1}^{n}w_{i}{\mbox
e}^{-\lambda(i-1)\Delta}](1-{\mbox e}^{-\lambda\Delta}){\mbox
e}^{y}{\mbox { for }}y\leq 0 \label{id1}\EndEq and for
$i=1,\ldots, n$\Eq
\sum_{j=1}^{n}w_{j}g_{j}{(y)}=\lambda^{-1}w_{i}(1-{\mbox e}^{-y})
{\mbox { for }}\lambda(i-1)\Delta<y\leq \lambda i\Delta.
\label{id2}\EndEq
 More generally suppose that for $t>s$,
$\tau(t)-\tau(s)=z(g_{s,t})=N(f_{s,t})$ where
$f_{s,t}(u,y)=ug_{s,t}(y)$ and \Eq
g_{s,t}(y)=h_{1,s,t,\lambda}(y)I_{\{\lambda s<y\leq \lambda
t\}}+h_{2,s,t,\lambda}F(y)I_{\{y\leq 0\}},\label{geng1}\EndEq for
$h_{1,s,t,\lambda}(y)$ and $F(y)$ non-negative functions
satisfying suitable integrability conditions and
$h_{2,s,t,\lambda}$ a positive quantity not depending on $y$.
Hence, one could choose for each $i$, \Eq
g_{i}(y)=h_{1,\lambda,i}(y)I_{\{\lambda(i-1)\Delta<y\leq \lambda
i\Delta\}}+h_{2,\lambda,i}F(y)I_{\{y\leq 0\}},\label{geng}\EndEq
for arbitrary positive functions $h_{1,\lambda,i}$, $F$ and
$h_{2,\lambda,i}$ whose form is determined by the general
difference $\tau(t)-\tau(s)$ for $t>s$. These models all exhibit
behavior similar to~\mref{id1} and ~\mref{id2}. That is for any
sequence of numbers $(w_{1},\ldots,w_{n})$, it follows that still
$
\sum_{i=1}^{n}w_{i}\tau_{i}=z(\sum_{i=1}^{n}w_{i}g_{i})=N(\sum_{i=1}^{n}w_{i}f_{i})$,
for $f_{i}(u,y)=ug_{i}(y)$ and $g_{i}$ now given by ~\mref{geng1}.
Additionally, the most important feature is preserved. That is,
\Eq \sum_{j=1}^{n}w_{j}g_{j}{(y)}=[\sum_{i=1}^{n}w_{i}
h_{2,\lambda,i}]F(y){\mbox { for }}y\leq 0 \label{id3}\EndEq and
for $i=1,\ldots, n$,
$\sum_{j=1}^{n}w_{j}g_{j}{(y)}=w_{i}h_{1,\lambda,i}(y) {\mbox {
for }}\lambda(i-1)\Delta<y\leq \lambda i\Delta.$

Now let $N$ denote a Poisson random measure on
$(0,\infty)\times(-\infty,\infty)$ with \Eq
\nu(dw,dy)=\rho_{1}(dw)\eta_{1}(dy)I_{\{y>0\}}+\rho_{2}(dw)\eta_{2}(dy)I_{\{y\leq
0\}},\label{nusim}\EndEq where $\rho_{1}$ and $\rho_{2}$ are
L\'evy densities generating L\'evy exponents $\psi_{1}$ and
$\psi_{2}$, and $\eta_{1},\eta_{2}$ are non-negative sigma-finite
measures. It follows from \mref{id3} that\Eq
\int_{0}^{\infty}\psi_{1}(\sum_{i=1}^{n}w_{i}g_{i}(y))\eta_{1}(dy)+
\int_{-\infty}^{0}\psi_{2}(\sum_{i=1}^{n}w_{i}g_{i}(y))\eta_{2}(dy)
=\Phi_{0}(s_{n})+\sum_{i=1}^{n}\Phi_{i}(w_{i})
\label{lapmix}\EndEq where $$
\Phi_{i}(w_{i})=\int_{\lambda(i-1)\Delta}^{i\Delta}
\psi_{1}(w_{i}h_{1,i,\lambda}(y))\eta_{1}(dy){\mbox { and
}}\Phi_{0}(s_{n})=\int_{-\infty}^{0}
\psi_{2}([\sum_{i=1}^{n}w_{i}h_{2,i,\lambda}]F(y))\eta_{2}(dy)
$$
for $i=1,\ldots,n$, where $s_{n}=
\sum_{i=1}^{n}w_{i}h_{2,i,\lambda}.$ In the case of~\mref{id1}
and~\mref{id2} for $\rho_{1}=\rho_{2}$, and
$\eta_{1}(dy)=\eta_{2}(dy)=dy$, $s_{n}=[\sum_{i=1}^{n}w_{i}{\mbox
e}^{-\lambda(i-1)\Delta}](1-{\mbox e}^{-\lambda \Delta})$. One has
for $i=1,\ldots,n$,\Eq \Phi_{i}(w_{i})=\int_{(1-{\mbox
e}^{-\lambda (i-1)\Delta})}^{(1-{\mbox e}^{-\lambda i\Delta})}
\psi(\lambda^{-1}w_{i}u) \frac{du}{1-u}=\int_{\lambda
(i-1)\Delta}^{\lambda i\Delta}\psi(w_{i}\lambda^{-1}(1-{\mbox
e}^{-y}))dy\label{levyv1}\EndEq and
$\Phi_{0}(s_{n})=\int_{0}^{1}\psi(\lambda^{-1}s_{n}u)\frac{du}{u}$
is the L\'evy exponent corresponding to the {\it prior
distribution} of $v_{0}$ evaluated at $s_{n}$.

\begin{thm}Let $N$ denote a Poisson random measure on $(0,\infty)\times(-\infty,\infty)$
with intensity $\nu$ defined in~\mref{nusim}. Define $\tau$
by~\mref{geng1} and ~\mref{geng} and the general specification of
$\nu$ above. Let $\Phi_{j}$ for $j=0,1,\ldots,n$ denote the
quantities defined by~\mref{lapmix}. Let $(v_{0},T_{1},\ldots,
T_{n})$ denote $n+1$ independent random variables with respective
Laplace transforms ${\mbox e}^{-\Phi_{i}(\omega)}$ for
$i=0,1,\ldots,n$.  Then, \Enumerate
\item[(i)]$\E\[{\mbox e}^{-\sum_{i=1}^{n}w_{i}\tau_{i}}\]={\mbox
e}^{-\Phi_{0}(s_{n})}\prod_{i=1}^{n}{\mbox
e}^{-\Phi_{i}(w_{i})}=\E[{\mbox
e}^{-s_{n}v(0)}|\nu]\prod_{i=1}^{n}\E\[{\mbox
e}^{-w_{i}T_{i}}|\nu\]$
\item[(ii)]
$\Lcr(\X|\mu,\beta,\theta)=
\int_{0}^{\infty}\[\prod_{i=1}^{n}\int_{0}^{\infty}
\phi(X_{i}|\mu\Delta+\beta z_{i},z_{i})
f_{T_{i}}(z_{i}-b_{i}y)dz_{i}\]f_{v_{0}}(y)dy$ where
$b_{i}=h_{2,i,\lambda}$ for $i=1,\ldots,n$.
$(f_{v_{0}},f_{T_{1}},\ldots,f_{T_{n}})$ denotes the densities of
the corresponding random variables.\qed \EndEnumerate
\end{thm}
\Proof  The result is a special case of Theorem 2.3 where one
replaces $\sum_{i=1}^{n}w_{i}\tau_{i}$ with
$\sum_{i=1}^{n}w_{i}[b_{i}v_{0}+T_{i}]$.\EndProof

 Bayesian inference may be conducted using the following
result.
\begin{prop} Suppose that $(T_{1},\ldots, T_{n},v_{0})$ depend on
external parameters, say $\theta$. Then assuming independent
priors $q(d\theta)$, $q(d\mu)$ and $q(d\beta)$, posterior
inference may be obtained based on the model derived from
augmenting the likelihood in Theorem 5.1. That is, the joint
distribution of $(\X,T_{1},\ldots,T_{n},v_{0},\theta,\mu,\beta)$
given by
$$
\[\prod_{i=1}^{n}\phi(X_{i}|\mu\Delta+\beta
Z_{i},Z_{i})\]q(d\mu)q(d\beta)
\[\prod_{i=1}^{n}f_{T_{i}}(Z_{i}-b_{i}v_{0})\]f_{v_{0}}(v_{0})q(d\theta).
$$
\qed
\end{prop}

\subsection{Predictive density of the stock price}
The model~\mref{BNS} suggests for any time interval $[s,t]$ that
the return of the log stock price, say $X_{s,t}$, is given by the
model, \Eq
X_{s,t}=\mu(t-s)+[\tau(t)-\tau(s)]\beta+{[\tau(t)-\tau(s)]}^{1/2}\epsilon_{s,t}
\label{genreturn} \EndEq where $\epsilon_{s,t}$ is an independent
standard Normal distribution. Our results yield an explicit
tractable expression of a predictive density of $X_{s,t}$ given
previously observed data $\X$

\begin{prop} For $t>s>n\Delta$, let $X_{s,t}$ be defined according to~\mref{genreturn}.
Let $\tau$ be defined by the general specifications in Theorem 5.1
Let $f_{T_{1}},\ldots,f_{T_{n}},f_{v_{0}}$ denote the densities of
the corresponding independent random variables. Then the
predictive density of $X_{s,t}|X_{1},\ldots, X_{n}$ is given by
the formula,
$$
\int_{0}^{\infty}\[\int_{0}^{\infty}\phi(x|\mu(t-s)+\beta
z_{s,t},z_{s,t}) f_{T_{s,t}}(z_{s,t}-b_{s,t}v)dz_{s,t}\]r(\X|v)
f_{v_{0}}(v)dv $$ where,
$r(\X|v)=\[\prod_{i=1}^{n}\int_{0}^{\infty}
\phi(X_{i}|\mu\Delta+\beta z_{i},z_{i})
f_{T_{i}}(z_{i}-b_{i}v)dz_{i}\]/\Lcr(\X|\mu,\beta,\theta)$ and
$b_{s,t}=h_{2,s,t,\lambda}$. The quantity, $f_{T_{s,t}}$ denotes
the density of an independent random variable $T_{s,t}$ with law
determined by the L\'evy exponent $\Phi_{s,t}(w)=\int_{\lambda
s}^{\lambda t} \psi_{1}(w h_{1,i,\lambda}(y))\eta_{1}(dy).$
\qed\end{prop} \Proof The result follows by noting that
$\sum_{i=1}^{n}w_{i}g_{i}(y)+w_{n+1}g_{s,t}(y)$ has the same
structural form as~\mref{id3}. Where $g_{s,t}$ is defined
in~\mref{geng1}. \EndProof
\subsection{Example}
The expressions in Theorem 4.1 suggests that an easily analyzed
model would arise if  $(v_{0},T_{1}\ldots,T_{n})$ were all from
GIG class of densities. Here, going back to the variant of the
Barndorff-Nielsen and Shephard model characterized
by~\mref{levyv1}, we shall show that the choice of a stable law
yields very nice results. Recall that the L\'evy exponent of a
stable law of index $0<\alpha<1$, is such that
$\psi(\omega)=\omega^{\alpha}/\alpha.$ Recall also that the case
of $\alpha=1/2$ leads to the Inverse Gamma distribution of index
$1/2$, and that the Inverse Gaussian arises from an exponential
tilting of this law. Using this fact we arrive at the following
result.
\begin{prop} Suppose that $\tau$ is specified by~\mref{model2}
with $\nu(ds,dy)=s^{-\alpha-1}/[\Gamma(1-\alpha)]dsdy$ for
$0<\alpha<1$. Then the random variables
$(v_{0},T_{1}\ldots,T_{n})$, appearing in Theorem 5.1, are
independent stable random variables of index $\alpha$. The
respective L\'evy exponents are, \Enumerate
\item[(i)]
$\Phi_{0}(\omega)=\omega^{\alpha}\lambda^{-\alpha}{(1-{\mbox
e}^{-\lambda \Delta})}^{\alpha}/\alpha $ and
\item[(ii)]
$\Phi_{i}(\omega)=\omega^{\alpha}\lambda^{-\alpha}\int_{1-{\mbox
e}^{-\lambda (i-1)\Delta}}^{1-{\mbox e}^{-\lambda
i\Delta}}\frac{u^{\alpha}}{\alpha(1-u)}du$ for $i=1,\ldots,n$
\item[(iii)] If one instead uses $\eta_{1}(dy)={\mbox
e}^{-y}dy$, then
$$\Phi_{i}(\omega)=\frac{\omega^{\alpha}}{\alpha(\alpha+1)}\lambda^{-\alpha}[{(1-{\mbox
e}^{-\lambda i\Delta})}^{\alpha+1}-{(1-{\mbox e}^{-\lambda (i-1)
\Delta})}^{\alpha+1}],$$ for $i=1,\ldots,n.$\EndEnumerate \qed
\end{prop}

Notice that the proposition above shows that a change from
Lebesque measure to $\eta_{{1}}(dy)$ only affects the constants in
the Laplace transform and generally preserves the distributional
property of the $(T_{i})$. This fits into what has been evidenced
in Bayesian nonparametric problems where the choice of quantities
such as $\eta_{1}$, $\eta_{2}$ are done mainly for computational
convenience. The rationale is that viewing the specifications for
$\tau$ as a prior model, experience from the Bayesian
nonparametric literature suggests that many such choices of $\tau$
will eventually lead to the similar conclusions in the presence of
enough data $\X$. Note however that we do not advocate removing
the dependence of the L\'evy exponents on $(i,\Delta)$, as this is
related to the data. We close by noting it is always possible to
arrange for the random variables $(v_{0},T_{1},\ldots,T_{n})$ to
be {\it self-decomposable} by choosing $\psi_{1}$ and $\psi_{2}$,
$\eta_{1}$ and $\eta_{2}$ such that the random variables are {\it
Generalized Gamma Convolutions}~(GGC). See Thorin~(1977) and
Bondesson~(1979, 1992) for this rich class of models. That is
$\Phi_{i}(\omega)=\int_{a_{i}}^{c_{i}}\ln(1+\omega/y)\Ucr(dy)$,
for $a_{i}$,$c_{i}$ depending on $\Delta$. $\Ucr$, with
$\Ucr(0)=0$, is called a Thorin function or measure. This applies
more generally to the models described in section 2.4. We shall
leave it to the reader to investigate which of the general models
discussed in section 2 are most suitable to their particular
application.\section {Extension: SV Likelihood models with
correlated jumps in price, leverage effects models.} Recall that
the Barndorff-Nielsen and Shephard(2001a,b) OU process $v(t)$,
which models the instantaneous volatility, satisfies the
differential equation $$dv(t)=-\lambda v(t)+dz(\lambda t),$$ where
the process $z$ is defined in section 3, and hence the volatility
possesses jumps. An important extension of the model in~\mref{BNS}
and hence to our general framework described in section 2, is
where one includes jumps in the log-price model which are
correlated with the the volatility $v$. These types of continuous
time models, wherein Duffie, Pan and Singleton~(2000) is an early
reference, serve to incorporate the {\it leverage effect }
discussed in for instance Black~(1976) and Nelson~(1991). We shall
be rather brief on this growing literature and refer the reader to
the works of Eraker, Johannes and Polson~(2003), Duan, Ritchken
and Sun~(2004) and Duffie, Singleton and Pan~(2000), for more
extensive background and rationale for these types of models and
its parametric variations. \subsection{BNS-OU SV model for
leverage effects} Barndorff-Nielsen and Shephard~(2001a, eq. 8)
describe this type of extension as follows,

\Eq dx^{*}(t)=(\mu+\beta v(t))dt+v^{1/2}(t)dw(t)+\rho (d[z(\lambda
t)-\E[dz(\lambda(t)]) \label{BNSLev}\EndEq assuming of course that
$\E[z(\lambda t)]<\infty$. One can incorporate modifications to
relax this condition. It follows that obviously the log price and
the volatility are negatively correlated if $\rho<0$. Thus
modeling the leverage effect that a fall in price results in an
increase in future volatility. Barndorff-Nielsen and
Shephard~(2001a, section 4) discuss further details of this model.
The likelihood model based on~\mref{BNSLev} was not explicitly
discussed in that paper, and moreover is considered even more
challenging. However as we shall show, this extension and a
variety of natural extensions of the models described in section
2, incorporating a leverage type effect,  are easily handled by
the type of methodology we have presented so far.

First, assuming a similar framework as in section 1.1, and using
the BNS-OU model described in section 3, note that
$$z_{i}:=z(g_{i,3})=z(i\Delta \lambda)-z(\lambda(i-1)\Delta)=\int_{0}^{\infty}\int_{0}^{\infty}I_{\{(i-1)\Delta\lambda<y\leq \lambda
i\Delta\}}uN(du,dy),$$ where $N$ is a Poisson random measure with
intensity $\nu(du,dy)=\rho(du)dy$ on $(0,\infty)\times
(-\infty,\infty)$, and $g_{i,3}(y)=I_{\{(i-1)\Delta\lambda<y\leq
\lambda i\Delta\}}.$ Assuming a finite first moment, one has
$\E[z_{i}]=\Delta\int_{0}^{\infty}u\rho(du)$. Hence the
model~\mref{BNSLev} implies that $X_{i}|\tau_{i},z_{i},\beta, \mu$
are conditionally independent with \Eq
X_{i}=\mu\Delta+\tau_{i}\beta+\tau^{1/2}_{i}\epsilon_{i}+\rho(z_{i}-\E[z_{i}]),
\label{dataJump}\EndEq which may be rewritten  for each $i$, as
$X_{i}=(\mu\Delta-\rho\E[z_{i}])+(\tau_{i}\beta+\rho
z_{i})+\tau^{1/2}_{i}\epsilon_{i}$. Hence on may write the
expression in~\mref{dataJump} as, \Eq
X_{i}=(\mu\Delta-\rho\E[z_{i}])+(z(g_{i,1}+g_{i,2})\beta+\rho
z(g_{i,3}))+\sqrt{z(g_{i,1}+g_{i,2})}\epsilon_{i},
\label{dataJump2}\EndEq which obviously may be further expressed
in terms of a common Poisson random measure.
\subsection{A General class of likelihoods which incorporate
leverage type effects} We note that from our point of view a model
such as ~\mref{dataJump2} poses no additional complications.
Similar to section 2, we will obtain exact  expressions for
likelihoods of quite general extensions of models with correlated
jumps in price and volatility. As before, for a general Poisson
random measure on $\Vcr$, with intensity $\nu$, let
$\tau_{i}=N(f_{i})$ for $i=1,\ldots,n$. Additionally, for
real-valued functions $\varphi_{1},\ldots,\varphi_{n}$ on $\Vcr$,
each satisfying the condition $\Lambda(|\varphi_{i}|)<\infty$, for
$i=1,\ldots,n,$ define, $\gamma_{i}:=N(\varphi_{i}).$ Now a
general version of ~\mref{dataJump2} is given by the case of
conditionally independent \Eq X_{i}=\mu\Delta+\rho
\gamma_{i}+\tau_{i}\beta+\tau^{1/2}_{i}\epsilon_{i}.
\label{dataJump3}\EndEq Let $\Scr={\mathbb R}\times{\mathbb
R}_{+}$, and assume that $N$ depends on a parameter $\theta.$ Then
the likelihood of $\X|\mu,\beta,\theta,\rho$ determined
by~\mref{dataJump3} can be expressed as \Eq
\Lcr(\X|\mu,\beta,\theta,\rho)=\int_{\Scr^{n}}
\[\prod_{i=1}^{n}\phi(X_{i}|\mu\Delta+\rho \gamma_{i}+\beta
\tau_{i},\tau_{i})\]f((\tau_{1},
\gamma_{1}),\ldots,(\tau_{n},\gamma_{n})|\theta)\prod_{i=1}^{n}d\gamma_{i}d\tau_{i}
\label{BNSliklev} \EndEq where,
$$\phi(X_{i}|\mu\Delta+\rho \gamma_{i}+\beta
\tau_{i},\tau_{i})={\mbox e}^{(A_{i}-\rho\gamma_{i})\beta
}\frac{1}{\sqrt{2\pi}}\tau^{-1/2}_{i}[{\mbox e}^{-{(A_{i}-\rho
\gamma_{i})}^{2}/(2\tau_{i})}]{\mbox e}^{-\tau_{i}\beta^{2}/2}
$$
The difficulty in evaluating this likelihood now manifests itself
in the term in brackets where both $\gamma_{i}$ and $\tau_{i}$ are
functionals of a common Poisson random measure, and moreover are
not pairwise independent across $i$. Clearly one could
apply~\mref{key1} however the cosine representation, now involving
random terms, would generally lead to expressions which are less
aesthetically pleasing. Here we will simply use an identity
deduced from the characteristic function of a Normal random
variable. This is equivalent to ~\mref{key1}. We close by
describing the explicit likelihood.
\begin{thm} Suppose that $N$ is a Poisson random measure with intensity $\nu$ on $\Vcr$. Then define
the complex valued function
$\Upsilon_{n}(x)=\sum_{i=1}^{n}[\rho(\beta+\xi
y_{i})]\varphi_{i}(x)$ on $\Vcr$, where $\xi$ denotes the
imaginary number. Then defining $\Omega_{n}(x)$ as in Theorem 2.1
with $w_{i}=y_{i}^{2}/2+\beta^{2}/2$ for $i=1,\ldots,n,$ the
likelihood~\mref{BNSliklev} can be expressed as
$$\Lcr(\X|\mu,\beta,\theta,\rho)= \frac{{\mbox e}^{n{\bar A}\beta
}}{{(2\pi)}^{n}}\int_{{\mathbb R}^{n}} {\mbox
e}^{-\Lambda(\Omega_{n}+\Upsilon_{n})} \prod_{i=1}^{n}{\mbox
e}^{\xi A_{i}y_{i}}dy_{i}.$$ where $(\varphi_{i})$ and $(f_{i})$
are chosen such that $\Lambda(\Omega_{n}+\Upsilon_{n})<\infty.$
 \qed\end{thm}
\Proof Here we use the fact that for each $i$ one has the identity
deduced from the characteristic function of a Normal distribution,
with mean 0 and variance $2/\tau_{i}$, evaluated at
$\varpi_{i}=A_{i}-\rho \gamma_{i}$. That is,
$$
\frac{1}{\sqrt{2\pi}}\int_{-\infty}^{\infty}{\mbox
e}^{\xi\varpi_{i}y_{i}-\tau_{i}y^{2}_{i}/2}dy_{i}=\frac{1}{\sqrt{\tau_{i}}}{\mbox
e}^{-{(\varpi_{i})}^{2}/2\tau_{i}}
$$
Now similar to the results in section 2, apply an appropriate
substitution, Fubini's theorem and the fact that $\rho
\gamma_{i}(\beta+\xi y_{i})+\tau_{i}w_{i}=N(\rho(\beta+\xi
y_{i})\varphi_{i}+w_{i}f_{i})$. That is after rearranging terms it
remains to calculate the expectation of ${\mbox
e}^{-N(\Omega_{n}+\Upsilon_{n})}$ \EndProof \Remark Again we note
that similar to the likelihood in Theorem 2.1, the likelihood
incorporating a generalized notion of leverage effects in Theorem
5.1 can be easily evaluated by classical numerical integration.
Additionally although we have concentrated on the
Barndorff-Nielsen and Shephard~(2001a, b) models, our framework
covers a large class of popular models in the literature, which
can now be be analyzed in a likelihood framework with leverage
effects. For some examples, see Carr, Geman,  Madan, and
Yor~(2003).We note further that since we used an identity that
does not depend on the distributional features of the Poisson
linear functionals the results can be easily adapted to other
processes with for instance possible additional Gaussian
components.\EndRemark

 \vskip0.2in {\Heading Acknowledgements}
This paper was heavily influenced by my interactions with John W.
Lau, to whom I extend my thanks. Thanks also to Sam Wong for his
support and stimulating conversation related to this topic. Thanks
to Albert Lo for pointing out the literature on models
incorporating leverage effects.

 \vskip0.2in \centerline{\Heading References}
\vskip0.2in \tenrm
\def\smc{\tensmc}
\def\sl{\tensl}
\def\bf{\tenbold}
\baselineskip0.15in

\Ref \by Andrews, G., Askey, R. and  Roy, R. \yr 1999 \book
Special functions. Encyclopedia of Mathematics and its
Applications, 71 \publ Cambridge University Press \publaddr
Cambridge \EndRef

\Ref \by Barndorff-Nielsen, O.E. and Shephard, N. \yr 2001a \paper
Ornstein-Uhlenbeck-based models and some of their uses in
financial economics \jour \JRSSB \vol 63 \pages 167-241 \EndRef
\Ref \by Barndorff-Nielsen, O.E. and Shephard, N. \yr 2001b \paper
Modelling by L\'evy processes for financial econometrics. In
L\'evy processes. Theory and applications. Edited by Ole E.
Barndorff-Nielsen, Thomas Mikosch and Sidney I. Resnick. p.
283-318. Birkh\"auser Boston, Inc., Boston, MA \EndRef

\Ref \by Barndorff-Nielsen, O. E. and Shephard, N.\yr 2003 \paper
Integrated OU processes and non-Gaussian OU-based stochastic
volatility models \jour Scand. J. Statist. \vol 30 \pages 277-295
\EndRef

\Ref \by Benth, F. E., Karlsen, K. H. and Reikvam, K. \yr 2003
\paper Merton's portfolio optimization problem in a Black and
Scholes market with non-Gaussian stochastic volatility of
Ornstein-Uhlenbeck type \jour Math. Finance \vol 13 \pages 215-244
\EndRef

\Ref \by Black, F. \yr 1976 \paper Studies of stock price
volatility changes \jour Proc. Bus. Econ. Statist. Sect. Am.
Statist. Assoc. \pages 177-181\EndRef

\Ref \by Black, F. and Scholes, M. \yr 1973 \paper The pricing of
options and corporate liabilities \jour J. Polit. Econ. \vol 81
\pages 637-654 \EndRef

\Ref \by Bondesson, L. \yr 1979 \paper A general result on
infinite divisibility \jour \AnnProb \vol 7 \pages 965-979 \EndRef

\Ref \by Bondesson, L. \yr 1992 \paper Generalized gamma
convolutions and related classes of distributions and densities.
Lecture Notes in Statistics, 76. Springer-Verlag, New York \EndRef

\Ref \by Carr, P., Geman, H., Madan, D.B. and Yor, M. \yr 2003
\paper Stochastic volatility for L\'evy processes \jour Math.
Finance \vol 13 \pages 345-382 \EndRef

\Ref \by Carr, P. and Schr\"oder, M. \yr 2004 \paper Bessel
processes, the integral of geometric Brownian motion, and Asian
options \jour Theor. Probab. Appl. \vol 48 \pages 400-425 \EndRef

\Ref \by Daley, D. J. and Vere-Jones, D. \yr 1988 \book An
introduction to the theory of point processes \publ
Springer-Verlag \publaddr New York \EndRef

\Ref \by Devroye, L. \yr 1986a \paper An automatic method for
generating random variates with a given characteristic function
\jour SIAM J. Appl. Math. \vol 46 \pages 698-719\EndRef

\Ref \by Devroye, L. \yr 1986b \book Nonuniform random variate
generation. \publ Springer-Verlag \publaddr New York\EndRef

\Ref \by Duan,  J.\yr 1995 \paper The GARCH option pricing model
\jour Math. Finance \vol 5 \pages 13-32 \EndRef

\Ref \by Duan, J., Ritchken, P. and Sun, Z. \yr 2004 \paper
Approximating GARCH-jump models, jump-diffusion processes and
option pricing. Manuscript available at
http://www.rotman.utoronto.ca/~jcduan/ \EndRef
 \Ref \by Duffie, D., Pan, J. and Singleton, K.,\yr 2000 \paper Transform Analysis and
Asset Pricing for Affine Jump Diffusions \jour Econometrica \vol
68 \pages 1343-1376 \EndRef

\Ref \by Eberlein, E. \yr 2001 \paper Application of generalized
hyperbolic L\'evy motions to finance. In L\'evy processes. Theory
and applications. Edited by Ole E. Barndorff-Nielsen, Thomas
Mikosch and Sidney I. Resnick. p. 319-336. Birkh\"auser Boston,
Inc., Boston, MA \EndRef

\Ref \by Engle, R. F.\yr 1982 \paper Autoregressive conditional
heteroscedasticity with estimates of the variance of United
Kingdom inflation \jour Econometrica \vol 50 \pages 987-1007
\EndRef

\Ref \by Eraker, B., Johannes, M., and Polson, N. \yr 2003 \paper
The impact of jumps in volatility and returns. \jour Journal of
Finance \vol 68 \pages 1269-1300 \EndRef

\Ref \yr 2004 \by Ishwaran,~H. and James,~L.~F. \paper
Computational methods for multiplicative intensity models using
       weighted gamma processes: proportional hazards, marked point
       processes and panel count data
\jour  \JASA \vol   99 \pages 175-190 \EndRef \Ref \by James, L.F.
\yr 2002 \paper Poisson process partition calculus with
applications to exchangeable models and Bayesian nonparametrics.
arXiv:math.PR/0205093 \EndRef

\Ref \by James, L.F. \yr 2005 \paper Bayesian Poisson process
partition calculus with an application to Bayesian L\'evy moving
averages. To appear in {\it Annals of Statistics}\\Available at
http://ihome.ust.hk/$\sim$lancelot/ \EndRef

\Ref \by Merton, R.~C. \yr 1973 \paper Theory of rational option
pricing \jour Bell J. Econ. Mgemt. Sci. \vol 4 \pages 141-183
\EndRef

\Ref \by Nelson, D. B. \yr 1991 \paper Conditional
heteroskedasticity in asset pricing: a new approach \jour
Econometrica \vol 59 \pages 347-370 \EndRef

\Ref \by Nicolato, E. and Venardos, E. \yr 2003 \paper Option
pricing in stochastic volatility models of the Ornstein-Uhlenbeck
type \jour Math. Finance \vol 13 \pages 445-466 \EndRef

\Ref \by Roberts, G. O., Papaspiliopoulos, O. and Dellaportas, P.
\yr 2004 \paper Bayesian inference for non-Gaussian
Ornstein-Uhlenbeck stochastic volatility processes \JRSSB \vol 66
\pages 369-393 \EndRef

\Ref \by Thorin, O.\yr 1977 \paper On the infinite divisibility of
the lognormal distribution \jour Scand. Actuar. J. \vol 3 \pages
121-148 \EndRef

\Ref \by Watson, G. N. \yr 1966 \book A treatise on the theory of
Bessel functions. Paperback Edition. Cambridge Mathematical
Library. Cambridge University Press, Cambridge \EndRef

\Ref \by Yor, M. \yr 1992 \paper On some exponential functionals
of Brownian motion \jour Adv. in Appl. Probab. \vol 24 \pages
509-531 \EndRef
\medskip
\smc

\Tabular{ll}

Lancelot F. James\\
The Hong Kong University of Science and Technology\\
Department of Information and Systems Management\\
Clear Water Bay, Kowloon\\
Hong Kong\\
\rm lancelot\at ust.hk\\
\EndTabular

\end{document}